\documentclass{article}    % Enable this line and disable the
\usepackage{graphicx}          % Include this line if your
\usepackage{epsfig}
\usepackage{graphicx}
\usepackage{makeidx}
\usepackage{multicol}
\usepackage{verbatim}
\usepackage{subfigure}
\usepackage{mathrsfs}
\usepackage{snapshot}

\usepackage{amsmath}

\usepackage{crop}
\crop
\usepackage{amssymb}%,minilisp,mycltxdefs}%,amsthm}%,verbatim}
\usepackage{graphicx}

\def\E{\mathbb{E}}

\def\ba{\begin{array}}
\def\ea{\end{array}}
\def\bi{\begin{itemize}}
\def\ei{\end{itemize}}

\DeclareMathOperator{\Tr}{Tr}

\newtheorem{definition}{Definition}[section]

\newtheorem{theorem}{Theorem}[section]
\newtheorem{lemma}{Lemma}[section]
\newtheorem{conjecture}{Conjecture}[section]

\begin{document}

%----
%\begin{frontmatter}
\begin{comment}
----
\title{Identification of Finite Dimensional Linear Systems Driven by L\'evy processes} % Title, preferably not more

%\thanks[footnoteinfo]{This paper was not presented at any IFAC
%meeting. Corresponding author M.~T.~Cicero. Tel. +XXXIX-VI-mmmxxi.
%Fax +XXXIX-VI-mmmxxv.}

\author[GL]{L\'aszl\'o Gerencs\'er}\ead{gerencser.laszlo@sztaki.mta.hu},    % Add the
\author[MM]{M\'at\'e M\'anfay}\ead{manfay@sztaki.mta.hu}            % e-mail address

\address[GL]{MTA SZTAKI}  % Please supply
\address[MM]{MTA SZTAKI, Central European University}             % full addresses
%\address[Baiae]{The White House, Baiae}        % here.

\begin{keyword}                           % Five to ten keywords,
linear systems, stochastic systems, L\'evy processes, system
identification, financial modelling              % chosen from the IFAC
\end{keyword}                             % keyword list or with the
                                          % help of the Automatica
                                          % keyword wizard
----
\end{comment}
\title{\LARGE \bf Stability of hybrid L\'evy systems}
\date{}
\author{L\'aszl\'o Gerencs\'er and M\'at\'e M\'anfay}

\maketitle

Continuous-time stochastic systems have attracted a lot of
attention recently, due to their wide-spread use in finance for
modelling price-dynamics. More recently models taking into accounts shocks
have been developed by assuming that the return process is an
infinitesimal L\'evy process. L\'evy processes are also used to model the traffic in a
telecommunication network. In this paper we focus
on a particular technical problem: stability of time-varying stochastic
systems driven or modulated by a L\'evy process with discrete time
interventions, such as parameter or state resetting. Such systems
will be called hybrid L\'evy systems. They are hybrid in the sense
that jumps both in the dynamics may occur. The
peculiarity of our systems is that the jump-times are defined by a
more or less arbitrary point process, but there exists an {\it
asymmetry} in the system dynamics.  The
novelty of our model relative to the theory of switching
stochastic systems is two-fold. First, we allow slow time
variation of the parameters, in a stochastic sense, without any
statistical pattern, in the spirit of the classical stability
result of Desoer, see \cite{desoer}. Secondly, we allow certain jumps
(resetting) in the system parameters almost without any a priori
condition.

%\subclass{subclass}
%\end{keyword}
%\pagestyle{myheadings}
%\thispagestyle{plain}
%\markboth{TEX PRODUCTION}{USING SIAM'S MM \LaTeX\ MACROS}

%\newpage

\section{Introduction}

Continuous-time stochastic systems have attracted a lot of
attention recently, due to their wide-spread use in finance for
modelling price-dynamics. A widely used model for continuous-time
returns has been, since the works of L. Bachelier, Gaussian white
noise with drift. More recently models taking into accounts shocks
have been developed by assuming that the return process is an
infinitesimal L\'evy process. For long term modelling a more
suitable model is a stochastic system with poles close to $1$
driven by a L\'evy process, see \cite{MM+GL_ECC_13}.

L\'evy processes are also used to model the traffic in a
telecommunication network. Other areas where L\'evy processes are used in modeling: robotics, mechanical systems, biology \cite{colaneri},\cite{saridis}. A further potential application is
modelling the shocks received by the wheel of a car due to the
irregularity of the road surface.

Stochastic processes driven or modulated by a L\'evy process will
be called a L\'evy system. Description of real data in terms of
L\'evy systems is far from being settled. In this paper we focus
on a particular technical problem that proved to be fundamental in
the statistical analysis of continuous-time stochastic systems
driven by Gaussian white noise, see \cite{hybrid_cikk}. Ultimately it
is hoped that this technical result may contribute to the
development of a continuous-time recursive maximum likelihood
method for finite dimensional linear stochastic L\'evy systems,
along the lines of \cite{GL+PROKAJ-RML-ECC09}.

The problem is the stability analysis of time-varying stochastic
systems driven or modulated by a L\'evy process with discrete time
interventions, such as parameter or state resetting. Such systems
will be called hybrid L\'evy systems. They are hybrid in the sense
that jumps both in the dynamics and the state may occur. The
peculiarity of our systems is that the jump-times are defined by a
more or less arbitrary point process, but there exists an {\it
asymmetry} in the system dynamics, inasmuch jumps can occur only
one-way, after a period of slow variation, namely, back to a fixed
point.

We note that the well-developed theory of switching stochastic
systems, see \cite{colaneri}, does not cover the problem that we consider. The
novelty of our model relative to the theory of switching
stochastic systems is two-fold. First, we allow slow time
variation of the parameters, in a stochastic sense, without any
statistical pattern, in the spirit of the classical stability
result of Desoer, see \cite{desoer}. Secondly, we allow certain jumps
(resetting) in the system parameters almost without any a priori
condition.

The structure of the paper is as follows: in Section II we develop
the basic technical tools, such as the geometric drift condition
and the associated Lyapunov-function method, for the analysis of
time-invariant L\'evy systems, and provide estimates for higher
order moments of the Lyapunov-function. In Section III we present
the simplest version of an extension of Desoer's theorem. In
Sections IV we prove a stability result under parameter resetting.
%%% Add Caines and Meyn, Colaneri, J. Lygeros,

A few basic notions related to L\'evy processes will be given in the Appendix. Throughout the paper we use the following notations: $|{\cdot}|$
stands for the Euclidean norm, $\| \cdot \|$ stands for the
induced matrix norm. $[X,Y]_t$ denotes the quadratic variation
of semi-martingales $X_t$ and $Y_t.$%, i.e., if $X$ is an \ito/ process in the form
%of $dX_t=a_td w_t+b_tdt$, then $\q{X}_t$ is a shorthand notation
%for $\int_0^t a_s^2ds$. The indicator of a set $A$ is denoted by
%$\I[A]$. Finally, $t\wedge s$ is the minimum of $t$ and $s$.

%{\new Section title:}

\section{Preliminaries}

Consider the time-invariant linear stochastic system
\begin{equation}
  \label{eq:STATE-EQN-T-INVARIANT}
  dX_t =AX_{t-} dt+B d L_t + C d W_t
\end{equation}
where $L_t$ is L\'evy process that has finite variation and has no continuous part and $W_t$ is a standard Wiener process. Assume that $A$ and $B$ are time-independent constant matrices. $X_t \in \mathbb{R}^{n}, L_t \in \mathbb{R}^{l}, W_t \in \mathbb{R}^{k}  A \in \mathbb{R}^{n \times n}, B \in \mathbb{R}^{n \times l}, C \in \mathbb{R}^{n \times k}.$  We will denote the $i^{th}$ component of a vector $V$ with $V^{(i)}.$ We will use the fact that if $L_t$ is of finite variation then it can be written as $L_t=L_0+bt+\sum_{s \leq t} \Delta L_s,$ and the quadratic covariance of two coordinates of such vector processes is of the form
$$\sum_{s \leq t} \Delta L^{(i)}_s \Delta L^{(j)}_s.$$
\begin{definition}
We say that a vector $L_t=(L_t^{(1)},\ldots,L^{(l)}_t)$ with independent components that are L\'evy processes satisfies the moment condition of order $Q$ if $\int_{\mathbb{R}}|x|^q \nu^{(i)}(dx)<\infty$ holds for all $1 \leq i \leq l,$ and for $1 \leq q \leq Q,$ where the L\'evy measure of $L^{(i)}_t$ is denoted by $\nu^{(i)}(x).$
\end{definition}

The next definition is motivated by the geometric drift condition introduced in $\cite{hybrid_cikk}.$ To estimate the moments of $X_t$ we will use a quadratic Lyapunov function $V_t.$

\begin{definition}
Let $L^{(i)}_t, 1 \leq i \leq l$ be independent L\'evy processes with finite variation. Let $f$ be a polynomial with coefficients bounded uniformly in $t$ and ${\rm deg} f \leq Q.$
Given a process $V_t$ satisfying with some $\varepsilon > 0$
\begin{equation}
\label{eq:drift_condition}
\begin{split}
dV_t=V_{t-} \Big(u_t dt + dM_t\Big) + \\
%\sum_{0 \leq j_1+\ldots+j_l \leq Q} c_{j_1,\ldots,j_l} \prod_{i=1}^{l} \left(\Delta L^{(i)}_t\right)^{j_i}\Big),
V_t^{1-\varepsilon}f(\Delta L^{(1)}_t,\ldots,\Delta L^{(l)}_t),
\end{split}
\end{equation}
where $\Delta L^{(i)}_t$ denotes the size of the jump of $L^{(i)}$ at $t.$ We say that $V_t$ satisfies the modified geometric drift condition of order $Q$ if there exist $\alpha,\gamma >0,$ such that
\begin{eqnarray}
\nonumber
  u_t &\leq& -\alpha \\
  \frac{d[M]_t}{dt} &\leq& \gamma.
\end{eqnarray}
\end{definition}

Without loss generality may always assume that decomposition of $L_t$ contains no drift term. Any possible drift term can be incorporated into $u_t dt.$ The following lemma will be used several times in the paper.
\begin{lemma}\label{P_lemma}
Assume Condition 1, let $D_0 \subset D$ compact and $\theta_0 \in D_0$. Then, there exists a smooth
function $P(\theta), \theta \in D$ and $\alpha > 0$ such that $P(\theta) \geq P(\theta_0) \geq I$ for all $\theta \in D_0$ and
$P(\theta)A(\theta) + A^T (\theta)P(\theta) \leq -\alpha P(\theta)$, for all $\theta \in D_0.$
\end{lemma}
For its proof see $\cite{hybrid_cikk}.$

The next two lemmas show that our Lyapunov function $V_t$ and its $q^{th}$ power satisfy the modified geometric drift condition.
\begin{lemma}\label{lemma:x^2}
Let $X_t$ be defined via $(\ref{eq:STATE-EQN-T-INVARIANT}),$ and $P$ given by Lemma $\ref{P_lemma}.$ Define $V_t=1+X^{T}_t P X_t,$ then $V_t$ satisfies the modified geometric drift condition of order two.
\end{lemma}

\begin{lemma}\label{V^q_cond}
%Assume that $L_t$ satisfies the moment condition of order $q$.
Let $X_t$ be defined via $(\ref{eq:STATE-EQN-T-INVARIANT}).$ Define $V_t=1+X^{T}_t P X_t,$ then $V^{q}_t$
satisfies the modified geometric drift condition of order $2q.$
\end{lemma}

The result of the next Lemma will be used in the proof of Theorem \ref{thm_^q_drift}.
\begin{lemma}
\label{lemma_^q_drift}
Let us suppose that $V_t$ satisfies the modified geometric drift condition of order $Q,$ and suppose that $L_t$ satisfies the moment condition of order $Q$ . Then
\begin{displaymath}
\E \left[V_t\right] < \infty
\end{displaymath}
holds.
\end{lemma}

{\bf Proof:}
Take a general process of the form $(\ref{eq:drift_condition}).$ Then $V_t$ satisfies
\begin{equation}\label{V_dolean}
\Delta V_t=V_{t-}Z_t,
\end{equation}
with
$$ Z_t=u_t dt+f(\Delta L^{(1)}_t,\ldots,\Delta L^{(l)}_t).$$
Using the Doleans-Dade exponential formula, see \cite{jacod}, for processes with finite variation yields the solution for $V_t:$
\begin{displaymath}
V_t = e^{ Z_t^{(c)} - Z_0^{(c)}} \prod_{s \leq t} \left(1 + \Delta Z_s \right),
\end{displaymath}
where $.^{(c)}$ denotes the continuous part of a process. This $V_t$ is also called as the stochastic exponential of $U_t.$
Let $c$ and $M$ be uniform bounds for $u_t$ and for the coefficients of $f.$ Increasing both $u_t$ and the coefficients of $f$  and taking absolute value of the jumps we obtain a bound on the solution $V_t:$
\begin{equation}\label{prod_est}
\begin{split}
V_t \leq e^{ct} \prod_{s \leq t} \left( 1 + M \sum_{0 \leq j_1+\ldots+j_l \leq Q} \prod_{i=1}^{l} \left|\Delta L^{(i)}_s\right|^{j_i} \right) \leq \\
e^{ct} \prod_{i=1}^l \prod_{s \leq t} \left( 1+M \sum_{j=1}^Q \left|\Delta L^{(i)}_s \right|^j  \right)
\end{split}
\end{equation}
Since $L^{(i)}_t$-s are independent processes it is sufficient show that
$$
\E\left[ \prod_{ s \leq t} \left(1+M \sum_{j=1}^Q \left|\Delta L^{(i)}_s \right|^j\right)\right]<\infty.
$$
Since $L_t$ satisfies the moment condition of order $Q$
$$
\E\left[1+M \sum_{j=1}^Q \left|\Delta L^{(i)}_s \right|^j\right] <\infty.
$$
Hence applying Lemma \ref{prod_exp}, see Appendix, concludes the proof.

The next theorem implies the stability of $X_t$ defined in $(\ref{eq:STATE-EQN-T-INVARIANT}).$
\begin{theorem}
\label{thm_^q_drift}
Let us suppose that $V_t$ satisfies the modified geometric drift condition of order $Q,$ and suppose that $L_t$ satisfies the moment condition of order $Q$ . Then
\begin{displaymath}
\sup_{t \geq 0}\E \left[V_t\right] < \infty
\end{displaymath}
holds for $1 \leq q \leq Q.$
\end{theorem}

The proof will be given in the Appendix.

\emph{Corollary:}
Since $V_t=1+X^{T}_t P X_t,$ satisfies the modified geometric drift condition so does $V^q_t.$ Hence, by Theorem \ref{thm_^q_drift} $\sup_{t \geq 0}\E \left[V^q_t\right] < \infty,$ which implies that $\sup_{t \geq 0}\E \left[|X|^q_t\right] < \infty$ holds under conditions seen in the theorem above.

\section{A stochastic Desoer's Theorem}

Consider a parametric family of linear stochastic state-space systems given by the state space equations:
\begin{equation}\label{X_dynamics}
dX_t=A(\theta_t)X_t dt+B(\theta_t)dW_t+C(\theta_t)dL_t
\end{equation}

Condition 1: $A(\theta)$ is stable for each $\theta \in D,$ where $D \subset \mathbb{R}^p$ is an open set, and $A(\theta),B(\theta)$ and $C(\theta)$ are smooth in $D$.

\begin{definition}
We say that $\theta_t$ is slowly varying in a stochastic sense if
\begin{equation}\label{slowly_var}
  d\theta_t=\beta_t dt + \sigma_t dW_t + \rho_t dL_t,
\end{equation}
with $|\beta|^2_t+||\sigma||^2_t+ ||\rho||^2_t < \delta,$ for some $\delta>0$ and all $t.$
\end{definition}

\begin{theorem}
Assume that $L_t$ satisfies the moment condition of order $Q$, $\theta_t$ is slowly
varying in the stochastic sense above, furthermore assume that $\theta_t \in \mathbb{R}^p$ is an adapted process taking its
values in a compact set $D_0 \subset D,$ and that Condition 1 holds. Then for a sufficiently small $\delta.$
\begin{equation}\label{desoer_claim}
  \sup_{t \geq 0} \E \left[ V_t^q\right] < \infty,
\end{equation}
for $1 \leq q \leq Q.$
\end{theorem}

{\bf Proof:}
The case when no L\'evy terms are present in the dynamics of $x_t$ and $\theta_t$ has been settled in Theorem 1 of  $\cite{hybrid_cikk}.$ We may therefore assume that $B(\theta)=0$ and $\sigma_t=0.$

For a given $\theta,$ let $P(\theta) \in C^2$ be a symmetric, positive definite matrix that solves
\begin{equation}\label{P_ineq}
  P(\theta)A(\theta)+A(\theta)^T P(\theta) \leq -\alpha P(\theta),
\end{equation}
with some $\alpha > 0,$ and $P(\theta)\geq I.$ Let $P_t=P(\theta_t),$ and consider $V_t=(1+X^T_t P_t X_t)^{q/2}$. It is enough to prove that $V_t$ satisfies the modified geometric drift condition. By Lemma $\ref{lemma_^q_drift}$ we only need to check that $1+X^T_t P_t X_t$ satisfies the modified geometric drift condition. We can write the dynamics of $1+X^T_t P_t X_t=1+\Tr(P_t Z_t),$ with $Z_t=X_t X^T_t$ as
\begin{equation}\label{varying dynamic}
\begin{split}
  d \Tr(P_t Z_t)=\Tr(P_t dZ_t)+\Tr(dP_t Z_t)+\\
  \sum_{i,j} dP_{i,j} d [X^{(i)},X^{(j)}]_t
\end{split}
\end{equation}
The first term can be handled using Lemma $\ref{lemma:x^2}.$ The dynamics of $P_t$ is given by
\begin{equation}\label{P_dynamics}
 dP_t=u_t dt + \Sigma_t d L_t,
\end{equation}
with $||u||^2_t+||\Sigma||^2_t < c\delta,$ with some $c.$ Thus, the second and the third term give drift terms that do not spoil the modified geometric drift condition.
The typical form of the contribution of the second term up to a bounded constant multiplier is
\begin{equation}\label{2._cont}
X^{(i)} X^{(j)} dL^{(k)},
\end{equation}
and that of the third term is
\begin{equation}\label{3._cont}
dL^{(k)}_t d[L^{(i)},L^{(j)}]_t.
\end{equation}
Hence, $1+X^T_t P_t X_t$ indeed satisfies the modified geometric drift condition. Thus, applying Theorem \ref{thm_^q_drift} concludes the proof.

This result implies the stability of the parameter varying system defined by ($\ref{X_dynamics}$).

\emph{Corollary:}
Under conditions seen in the previous theorem $\sup_{t \geq 0} \E \left[ |X|_t^q\right] < \infty$ holds.

\section{Jumps in the dynamics of the parameter}

We now assume that the slowly parameter varying process $\theta_t$ resets at random times defined by a point process with counting process $N_t.$
\begin{equation}\label{theta_resetting}
  d\theta_t=\beta_tdt+\sigma_t dW_t + dL_t+(\theta_0-\theta_{t-})dN_t,
\end{equation}
where $|\beta_t|^2+||\sigma_t||^2 < \delta.$
\begin{theorem}
Assume that Condition 1 holds, and that $L_t$ satisfies the moment condition of order $Q,$ and let $X_t$ be defined via $(\ref{X_dynamics}).$ Then
$$ \sup_{t \geq 0} \E\left[ |X|_t^q \right]< \infty$$
holds for $1 \leq q \leq Q.$
\end{theorem}

{\bf Proof:}
We may assume that there is no diffusion part in the dynamics of $X_t$ and $\theta_t.$
Let $P(\theta)$ be defined by Lemma \ref{P_lemma} so that it attains its minimum on $D$ in $\theta_0.$ Define $V_t=(1+X_t^T P(\theta_t)X_t)^{q/2}.$ Let $\xi_t$ be the size of the jump at $t$ induced by the jump of $\theta,$ ie.
\begin{equation}
\xi_t=(1+X_t^T P(\theta_0)X_t)^{q/2}-(1+X_t^T P(\theta_t)X_t)^{q/2},
\end{equation}
using this notation the dynamics of $V_t$ can be written as
\begin{equation}\label{V_dyn_theta_jumps}
%dV_t=V_{t-} \left(u_t(V_t) dt + \sum_{0 \leq j_1+\ldots+j_l \leq Q} c_{j_1,\ldots,j_l}(V_t) \prod_{i=1}^{l} \left(\Delta L^{(i)}_t\right)^{j_i}\right)+v_t dt+\xi_t dN_t.
dV_t=V_{t-} U_t+\xi_t dN_t,
\end{equation}
with
$$
U_t=u_t dt + \sum_{0 \leq j_1+\ldots+j_l \leq Q} c_{j_1,\ldots,j_l} \prod_{i=1}^{l} \left(\Delta L^{(i)}_t\right)^{j_i}.
$$
By the minimality of $P(\theta_0),$ the jump term in $(\ref{V_dyn_theta_jumps})$ causes a non-positive jump in $V_t.$ Let $\psi_t$ be the stochastic exponential of $U_t,$
then
\begin{equation}
V_t=\psi_t V_0+ \int_{0}^{t} \psi_t \psi_s^{-1} \xi_s dN_s \leq \psi_t V_0.
\end{equation}
Since $\E \left[ \psi_t V_0 \right]<\infty$ is implied by Theorem $\ref{thm_^q_drift},$ we conclude the proof.

\section{Discussion: State resetting for jump processes}

Consider the hybrid linear system with jumps
\begin{equation}\label{state_resetting_eq}
  dX_t=AX_t+BdW_t+CdL_t+(X_0-X_{t-})dN_t,
\end{equation}
where $W_t$ is a Wiener process, and $N_t$ is a counting process.

\begin{conjecture}\label{state_re_thm}
 Suppose that $L_t$ satisfies the moment condition of order $Q$, then for $X_t$ defined by $(\ref{state_resetting_eq})$
 $$ \sup_{t \geq 0}\E[|X|_t^q] < \infty,$$
 holds for $1 \leq q \leq Q.$
\end{conjecture}

Our future work will focus on proving this conjecture. Although the main ideas of the proof are established some technical issues are still to be proven.

\section{Appendix}
L\'evy processes with finite variation are formally obtained via
\begin{equation}
\label{eq:Levy_def}
Z_t = \int_0^t \int_{{\mathbf R}^1} x N(ds, dx),
\end{equation}
where $N(dt,dx)$ is a time-homogeneous, space-time Poisson point
process, counting the number of jumps of size $x$ at time $t$. L\'evy processes are characterized by their L\'evy measures that can be defined via the intensity of $N(dt,dx):$
$$\E [N(dt,dx)] = dt \cdot \nu(dx),$$
where $\nu(dx)$ is the L\'evy-measure.
The quadratic variation of semimartingales $X$ and $Y$ is defined by the following process:
$$
d[X,Y]_t=d(XY)_t-X_{t-}{dY_t}-Y_{t-}dX_t.
$$
If $X=Y$ we get the quadratic variation of $X.$

Proof of Lemma \ref{lemma:x^2}:

{\bf Proof:}
This lemma is an extension of Lemma 8 in \cite{hybrid_cikk}, where no L\'evy processes are present in defining the dynamics of $V_t.$ Thus, for the sake of simplicity, we may omit the martingale $M_t$ from (\ref{eq:drift_condition}).

%For simplicity, we omit the martingale $M_t$ from (\ref{eq:drift_condition}), in the case of the presence of a martingale the following proof can be simply extended with the reasoning seen in Lemma 8 in \cite{hybrid_cikk}.

%Define $V_t=1+X^{T}_t P X_t=1+\Tr(PZ_t),$ where $Z_t=X_t X^T_t.$
Write $V_t=1+X^{T}_t P X_t=1+\Tr(PZ_t),$ where $Z_t=X_t X^T_t.$ The dynamic of $Z_t$ can be written as
\begin{equation}
\label{eq:Z_dinamics}
dZ_t=X_{t-}dX^{T}_t+dX_t X^{T}_{t-}+ B d[L,L]_t B^T,
\end{equation}
where $d[L,L]_t$ is an $l\times l$ matrix with entries representing quadratic covariances, that is $d[L,L]^{(i,j)}_t=d[L^{(i)},L^{(j)}].$ Equation $(\ref{eq:Z_dinamics})$ reads as
\begin{equation}
\begin{split}
X_{t-}\left( X^{T}_{t-}A^T dt+dL^T_t B^T\right)+\left(A X_{t-} dt+B dL_t\right)X^{T}_t \\ +B d[L,L]_t B^T
\end{split}
\end{equation}
Thus the dynamics of $V_t=1+\Tr(PZ_t)$ can be written as
\begin{equation}
\begin{split}
\left( PX_{t-}X^T_{t-}A^T+PAX_{t-}X^{T}_{t-}\right)dt+ \\ + PX_{t-} dL_t B^T+PB dL_t X^T_{t-}+PB d[L,L]_t B^T.
\end{split}
\end{equation}
So the $dt$ terms in the dynamics of $V_t=1+\Tr(PZ_t)$ we can write
\begin{equation}
\begin{split}
\Tr\left( X^T_{t-} A^T PX_{t-}+X^T_{t-}PAX_{t-}\right) \\ \leq -\alpha X^T_{t-} P X_{t-}= -\alpha V_{t-}+ \alpha,
\end{split}
\end{equation}
for the terms having $dL_t$
\begin{equation}
\begin{split}
&\Tr\left( P(X_{t-} dL_t B^T+B dL_t X^T_{t-})\right)= \\
&\Tr\left((P+P^T)B dL_t X^T_{t-}\right)= \\
&2X^T_{t-} (P+P^T)B dL_t=\psi^T_t dL_t,
\end{split}
\end{equation}
with $|\psi_t|^2=4X^T_{t-}(P+P^T)BB^T(P+P^T)X_{t-} \leq 4KV_{t-},$ with some fixed $K.$
Finally for the term with $d[L,L]_t$
\begin{equation}
\begin{split}
\Tr\left( PBd[L,L]_t B^T\right)=\sum^{l}_{i,j=1} c_{i,j} d[L^{(i)},L^{(j)}]_t= \\
\sum^{l}_{i,j=1} c_{i,j} \Delta L^{(i)}_t \Delta L^{(j)}_t,
\end{split}
\end{equation}
with some $c_{i,j}, 1 \leq i,j \leq l$ constants. It follows that the dynamic of $V_t$ can be written as
\begin{equation}
\begin{split}
dV_t=V_{t-}u_t dt+ \\
V^{1/2}_{t-}\left(\sum^{l}_{i=1}\frac{\psi^{(i)}_t}{V^{1/2}_{t-}} \Delta L^{(i)}_t+\sum^{l}_{i,j=1} \frac{c_{i,j}}{V^{1/2}_{t-}} \Delta L^{(i)}_t \Delta L^{(j)}_t\right),
\end{split}
\end{equation}
with uniformly bounded $u_t, \frac{\psi^{(i)}_t}{V^{1/2}_{t-}}, \frac{c_{i,j}}{V^{1/2}_{t-}}$  for any $1 \leq i,j \leq l.$

Proof of Lemma \ref{V^q_cond}:

{\bf Proof:}
The dynamics of $V^q_t$ can be written as
\begin{equation}
\begin{split}
dV^q_t=qV^{q-1}_{t-} dV_{t,(c)}+V^{q}_t-V^{q}_{t-}= \\
qV^{q-1}_{t-} u_t dt+\left(V_{t-}+\Delta V_t\right)^q-V^{q}_{t-}= \\
qV^{q}_{t-} u_t/V_{t-1} dt+ \sum_{k=1}^{q} {q \choose k} \left(\Delta V_t\right)^{k}V^{q-k}_{t-} \label{rhs},
%V_{t-}^q \left( \frac{q u_t}{V_{t-}} dt + \sum_{k=0}^{q-1} {q \choose k} V^{k-q}_{t-}\left(\Delta V_t\right)^{q-k} \right)
\end{split}
\end{equation}
with $u_t/V_{t-1}<\alpha.$ Using that
$$
\Delta V_t=V^{1-\varepsilon}_{t-}f(\Delta L^{(1)}_t,\ldots,\Delta L^{(l)}_t),
$$
we obtain that a typical jump term in (\ref{rhs}) reads as up to constant multiplier
%\begin{equation}
%\begin{split}
%qV^{q-1}_{t-} u_t dt+V^q_{t-}\sum_{k=0}^{q-1} {q \choose k} \left(\sum^{l}_{i=1}\psi^{(i)}_t \Delta L^{(i)}_t+\sum^{l}_{i,j=1} c_{i,j} \Delta L^{(i)}_t \Delta L^{(i)}_t\right)^q= \\
%V^q_{t-}\left( \bar{u}_t dt +\sum_{k=0}^{q-1} {q \choose k} \left(\sum^{l}_{i=1}\psi^{(i)}_t \Delta L^{(i)}_t+\sum^{l}_{i,j=1} c_{i,j} \Delta L^{(i)}_t \Delta L^{(i)}_t\right)^q\right),
%qV^{q-1}_{t-} u_t dt+V^q_{t-}\sum_{k=0}^{q-1} {q \choose k} U_t^q=
%V^q_{t-}\left( \bar{u}_t dt +\sum_{k=0}^{q-1} {q \choose k} U_t^q\right),
%\end{split}
%\end{equation}
\begin{equation}
V^{q-k \varepsilon}_{t-}f(\Delta L^{(1)}_t,\ldots,\Delta L^{(l)}_t)^k
\end{equation}
This implies that $V_t^q$ satisfies the modified geometric drift condition of order $2q.$

The next two technical Lemmas will be used in the proof of Theorem \ref{thm_^q_drift}.

\begin{lemma}\label{prod_exp}
Let $L_t$ be a L\'evy process with L\'evy measure $\nu.$ Suppose that a function $f$ satisfies
$$
\int_{\mathbb{R}} f(x) \nu(dx) <\infty,
$$
then
$$
\E\left[ \prod_{s \leq t} (1+f(\Delta L_s))\right] =e^{t\int_{\mathbb{R}} f(x)\nu(dx)}
$$
holds for any $t$.
\end{lemma}

{\bf Proof:}
First suppose that $L_t$ is a compound Poisson process with rate $\lambda,$ then the expected value of
$$\psi_t=\prod_{s \leq t} f(\Delta L_s)$$
can be estimated by conditioning on the number of jumps of $L_t.$ Let $N_t, J_t$ denote the number of jumps of $L_t$ on $[0,t],$ and the set of time indices when $L$ jumps on $[0,t].$ Define $D_t^n=\{ (t_1,\ldots,t_n):0\leq t_i \leq t, \text{ for all } 1\leq i \leq n\}.$
\begin{displaymath}
\begin{split}
\label{eq:psi_estimate}
&\E \left[\Psi_t\right]=\sum_{n=0}^{\infty} \E \left[\Psi_t | N_t= n \right] \mathbb{P} \left(N_t=n \right) =\\
 &\int_{D_t^n} \E[\Psi_t | N_t=n, J_t=\{t_1,\ldots,t_n\} ] \\ &P(N_t=N)dt_1\ldots dt_n= \\
  &\sum_{n=0}^{\infty} (m+1)^n e^{-\lambda t} \frac{\left(\lambda t \right)^{n}}{n!}=e^{ \lambda m t},
\end{split}
\end{displaymath}
where $m=\E \left[f(\Delta L_t)| L \text{ jumps at } t\right],$ and $P$ is the joint probability density of the jump times. For the general case define the truncated L\'evy measure
\begin{displaymath}
\nu^{\varepsilon}(x)=\frac{\textbf{1} \left[\left|x\right|>\varepsilon \right]\nu(x) }{\int_{\left|x\right|>\varepsilon}\nu (dx)},
\end{displaymath}
and let $L_t^{\varepsilon}$ be the L\'evy process with L\'evy measure $\nu^{\varepsilon}$. Then $L_t$ is the weak limit of $L_t^{\varepsilon}$ as $\varepsilon$ tends to zero.
\begin{eqnarray}
&& m_f^{\varepsilon}=\int f(x) \nu^{\varepsilon} (dx) \\
&& m^{\varepsilon}=\E \left[f(\Delta L^{\varepsilon}_t)| L^{\varepsilon} \text{ jumps at } t\right] \\
&& \lambda^{\varepsilon}=\int_{\left|x\right|>\varepsilon}\nu (dx)
\end{eqnarray}
writing (\ref{eq:psi_estimate}) for $L^{\varepsilon}$ yields
\begin{displaymath}
\E\left[ \prod_{s \leq t} (1+f(\Delta L^{\varepsilon}_s))\right] = e^{\lambda^{\varepsilon} (m^{\varepsilon}-1) t}
\end{displaymath}
Note that $\lambda^{\varepsilon} m_k^{\varepsilon}= \int_{\left|x\right|>\varepsilon} f(x)  \nu (dx)$, it follows that $e^{\lambda^{\varepsilon} m^{\varepsilon}t}$ has finite limit as $\varepsilon \rightarrow 0^{+}$ provided $\int_{\mathbb{R}} f(x) \nu(dx) <\infty$ which is the case. Hence,
$\E \left[ \psi_t \right]=e^{t\int_{\mathbb{R}} f(x)\nu(dx)}$ follows.

\begin{lemma}\label{int_cauchy_est}
Let the one dimensional process $L_t$ with L\'evy measure $\nu$ satisfy the moment condition of order $Q.$ Let $f$ be a polynomial with ${\rm deg} f\leq Q,$ and $f(0)=0.$ Then
$$
\E\left[ \int_0^t e^{-\alpha(t-s)} f(\Delta L_s)\right]=\frac{1-e^{-\alpha t}}{\alpha}\int_{\mathbb{R}} f(x) \nu(dx)
$$
\end{lemma}

{\bf Proof:}
First consider the case when $L_t$ is a compound Poisson process with intensity $\lambda$. Let $N_t, J_t$ denote the number of jumps of $L_t$ on $[0,t],$ and the set of time indices when $L$ jumps on $[0,t].$ Define $D_t^n=\{ (t_1,\ldots,t_n):0\leq t_i \leq t, \text{ for all } 1\leq i \leq n\}.$
\begin{align}
%\begin{split}
&\E\left[ \int_0^t e^{-\alpha(t-s)} f(\Delta L_s)\right]= \nonumber \\
&\sum_{n=0}^{\infty}\E\left[ \int_0^t e^{-\alpha(t-s)} f(\Delta L_s)| N_t=n\right]P(N_t=n). \label{eq_01}
%\end{split}
\end{align}
Calculating one term in the sum above
\begin{equation}
\begin{split}
&\int_{D_t^n}\E\left[ \int_0^t e^{-\alpha(t-s)} f(\Delta L_s) | N_t=n,J_t=\{t_1,\ldots,t_n\}\right]\\
&P(N_t=n)\frac{dt_1 \ldots dt_n}{t^n}= \\
&\int_{D_t^n}\sum_{i=1}^{n} e^{-\alpha(t-t_i)}\E\left[f(\Delta L_s)|N_t=n,J_t=\{t_1,\ldots,t_n\}\right]\\ &e^{-\lambda t} \frac{(\lambda t)^n}{n!}\frac{dt_1 \ldots dt_n}{t^n}=\\
&n\int_{D_t^{n-1}}\int_0^{t} e^{-\alpha(t-t_n)} \E\left[f(\Delta L_s)|t_1 \in J_t\right] e^{-\lambda t} \frac{(\lambda t)^n}{n!} \\
& \frac{dt_n}{t} \frac{dt_1 \ldots dt_{n-1}}{t^{n-1}}= \\
&ne^{-\lambda t} \frac{(\lambda t)^n}{n!} \E\left[f(\Delta L_s)|t_1 \in J_t\right]\frac{1-e^{-\alpha t}}{\alpha t}.
\end{split}
\end{equation}
Now using this result in (\ref{eq_01}) yields
\begin{equation}\label{f_exp}
\begin{split}
&\frac{1-e^{-\alpha t}}{\alpha t} \E\left[f(\Delta L_s)|t_1 \in J_t\right] \sum_{n=0}^{\infty} n e^{-\lambda t} \frac{(\lambda t)^n}{n!}= \\
& \frac{1-e^{-\alpha t}}{\alpha } \lambda \E\left[f(\Delta L_s)|t_1 \in J_t\right]
\end{split}
\end{equation}
For the general case define like in the proof of Lemma \ref{prod_exp} process $L^{\varepsilon}_t$ and its L\'evy measure $\nu^{\varepsilon}(dx),$ and $m^{\varepsilon}=\E\left[f(\Delta L_s)|t_1 \in J_t\right].$ Writing (\ref{f_exp}) for $L^{\varepsilon}_t$ we obtain
\begin{equation}
\begin{split}
&\E\left[ \int_0^t e^{-\alpha(t-s)} f(\Delta L^{\varepsilon}_s)\right]= \\
&\frac{1-e^{-\alpha t}}{\alpha } \lambda^{\varepsilon} \E\left[f(\Delta L^{\varepsilon}_s)|t_1 \in J_t\right]= \\
&\frac{1-e^{-\alpha t}}{\alpha } \int_{|x|>\varepsilon} f(x) \nu(dx).
\end{split}
\end{equation}
Since $L_t$ is the weak limit of $L^{\varepsilon}_t$ as $\varepsilon$ tends to zero, allowing $\varepsilon \rightarrow 0^+$ concludes the proof.

Proof of Theorem \ref{thm_^q_drift}:

{\bf Proof:}
Let $V_t$ satisfy the modified geometric drift condition
\begin{equation}
dV_t=u_t dt+V_t^{1-\varepsilon}f(\Delta L^{(1)}_t,\ldots,\Delta L^{(l)}_t).
\end{equation}
In the presence of a martingale in the dynamics of $V_t$ one can apply methods shown in \cite{hybrid_cikk}.
Applying Cauchy formula gives
\begin{equation}
V_t=\int_0^t e^{-\alpha(t-s)}V_s^{1-\varepsilon} f(\Delta L^{(1)}_t,\ldots,\Delta L^{(l)}_t)
\end{equation}
Now we estimate the expected values of $V_t$ using $V^*_t=\sup_{0\leq s\leq T} \E[V_s]$
\begin{align}
&\E[V_t] = \int_0^t \E \left[ e^{-\alpha(t-s)} V_t^{1-\varepsilon} f(\Delta L^{(1)}_t,\ldots,\Delta L^{(l)}_t) \label{0_1} \right] = \\
& \int_0^t e^{-\alpha(t-s)} \E[V_s^{1-\varepsilon}] \E\left[f(\Delta L^{(1)}_t,\ldots,\Delta L^{(l)}_t) \right] \leq \label{0_2} \\
& \int_0^t e^{-\alpha(t-s)} \E[V_s]^{1-\varepsilon} \E\left[f(\Delta L^{(1)}_t,\ldots,\Delta L^{(l)}_t) \right] \leq \label{0_3} \\
& (V_T^*)^{1-\varepsilon} \int_0^t e^{-\alpha(t-s)} \E\left[f(\Delta L^{(1)}_t,\ldots,\Delta L^{(l)}_t) \right]
\label{0_4} \leq \\
&(V_T^*)^{1-\varepsilon} \prod_{i=1}^l \int_0^t e^{-\alpha(t-s)} \E\left[g_i(\Delta L^{(i)})\right] \label{0_5}
\end{align}
In (\ref{0_1}) we used Fubini's theorem and the independency of $V_{s-}$ and $\Delta L^{(i)}_s.$ In (\ref{0_2}) H\"older inequality was applied. Finally, in (\ref{0_5}) we estimate $f$ with products of $g_{i}$ polynomials as we did in (\ref{prod_est}), clearly ${\rm deg } g_i \leq Q$ holds for all $i.$ Applying Lemma \ref{int_cauchy_est} gives for (\ref{0_5})
$$
(V^*_T)^\varepsilon \leq \prod_{i=1}^l \frac{1-e^{-\alpha t}}{\alpha}\int_{\mathbb{R}} g_i(x) \nu_i(dx),
$$
where $\nu_i$ is the L\'evy measure of $L^{(i)}_t.$ Since $1-e^{-\alpha t}<1$ we obtained a bound on $V^*_T$ that do not depend on $T.$ Hence, $\sup_{0 \leq t} \E[V_t] <  \infty,$ which concludes the proof.

% The following packages can be found on http:\\www.ctan.org
%\usepackage{graphics} % for pdf, bitmapped graphics files
%\usepackage{epsfig} % for postscript graphics files
%\usepackage{mathptmx} % assumes new font selection scheme installed
%\usepackage{times} % assumes new font selection scheme installed
%\usepackage{amsmath} % assumes amsmath package installed
%\usepackage{amssymb}  % assumes amsmath package installed

\end{document}